\documentclass[12pt]{article}
\usepackage{amsmath,fullpage,amsthm,amssymb,xcolor,graphicx}

\newtheorem{thm}{Theorem}[section]

\newtheorem{definition}{Definition}[section]
\newtheorem{remark}{Remark}[section]

\title{On spectrally optimal duals for $r$-erasures of frames generated by graphs}
\author{Deepshikha \thanks{Department of Mathematics, Shyampur Siddheswari Mahavidyalaya, University of Calcutta, West Bengal 711312, India. Email: dpmmehra@gmail.com }\  \and Aniruddha Samanta \thanks{Theoretical Statistics and Mathematics Unit, Indian Statistical Institute, Kolkata-700108, India. Email: aniruddha.sam@gmail.com} 
}
\date{\today}
\begin{document}
	\maketitle
	\baselineskip=0.25in
	\begin{abstract}
In \cite{DS3}, authors studied spectrally optimal dual frames for $1$-erasure and $2$-erasures of frames generated by graph. In this paper, we study spectrally optimal dual frames for $r$-erasures. We show that the spectral radius of the error operator of unitary equivalent frames is same with respect to their respective canonical dual frames.  We prove that if a frame is generated by a connected graph, then its canonical dual frame is the unique spectrally optimal dual frame for $r$-erasures. Further, we show that the canonical dual of frames generated by disconnected graphs are non-unique spectrally optimal dual frames for $r$-erasures.
	\end{abstract}

	
	
	\section{Introduction}
 In order to solve certain complex problems in the non-harmonic Fourier series, Duffin and Schaeffer \cite{DS} developed the idea of frames. Frames are redundant sequences of vectors that span a vector space. The redundancy of frames makes them very useful in many applications, including data and image processing, signal transmission, image segmentation, coding theory, etc. Due to redundancy, the data transmitted using frames can be reconstructed with minimal error if erasures occur during the transmission.  

 Finite frames are widely used in frame theory due to their significant relevance in applications. In \cite{D}, authors introduced a new class of finite frames using simple graphs and called them frames generated by graphs. They established some graph properties using frame properties and vice-versa. Consequently, it provides a beautiful connection between the three classical branches of mathematics, namely, frame theory, graph theory, and matrix theory. Further, frames generated by graphs have some special properties, and some of them are stated in the preliminaries section of this article.

 The problem of erasures commonly occurs due to transmission losses or disturbances during transmission. The problem of erasures can be handled in two ways. One way is to find a frame that can minimize the maximal error in reconstruction. Another way is to find an optimal dual frame for a pre-selected frame in order to minimize the maximal error of reconstruction. In this paper, we focus on the second approach. A dual frame of a pre-selected frame that minimizes the maximal error of reconstruction is called an optimal dual frame. Various measurements are used to measure the error operator. If operator norm is used, then we use the term optimal dual frame. If the spectral radius is used to measure the error operator, then dual frames that minimize the maximum error are called spectrally optimal dual frames (or, in short, $SOD$-frames). In \cite{DS1}, authors studied optimal dual frames of frames generated by graphs. Spectrally optimal dual frames for $1$-erasure and 
$2$-erasures of frames generated by graphs is studied in \cite{DS3}. In this paper, we study the spectrally optimal dualS for $r$erasures of frames generated by graphs. 
		\vspace{4pt}
		
	The outline of this paper is as follows. In Section 2, we provide the basic definitions and results of frame theory, spectrally optimal dual frames, spectral graph theory, frames generated by graphs and $SOD$-frame for $1$-erasure and $2$-erasures required in the rest of the article. In Section 3, we study spectrally optimal dual frames for $r$-erasures of frames generated by graphs.We show that the spectral radius of the error operator of unitary equivalent frames is same with respect to their respective canonical dual frames.  We prove that if $\Phi=\{\phi_i\}_{i\in[n]}$ is an $L_G(n,n-k)$-frame with frame operator $S_{\Phi}$, then $\rho^{(r)}_{\Phi, S^{-1}_{\Phi}\Phi}=1$ for any $r\in[n]\setminus\{1\}$. We prove that if a frame is generated by a connected graph, then its canonical dual frame is the unique spectrally optimal dual frame for $r$-erasures. Further, we show that the canonical dual of frames generated by disconnected graphs are non-unique spectrally optimal dual frames for $r$-erasures..
	
	\section{Preliminaries}
Throughout the rest of the article, we consider finite frames and $\mathbb{C}^k$ is a $k$-dimensional Hilbert space. We use the notation $[n]$ to denote the set $\{1,2,\ldots,n\}$ and $M=[m_{i,j}]_{i\in [p], j\in [q]}$ to denote the matrix $M$ of order $p\times q$ where $m_{i,j}$ is the $(i,j)$-th entry of $M$. Also, $diag(\alpha_1,\alpha_2,\ldots,\alpha_n)$ denotes the diagonal matrix of order $n$ with respective diagonal entries $\alpha_1,\alpha_2,\ldots,\alpha_n$. 

  A sequence of vectors $\Phi:=\{\phi_i\}_{i\in[n]}$ in $\mathbb{C}^k$ is called a \emph{frame} for $\mathbb{C}^k$ if there exist constants $0<A_o\leq B_o$ such that for any $f\in\mathbb{C}^k$, we have
\begin{align*}
A_o\|f\|^2\leq\sum_{i\in[n]}|\langle f, \phi_i\rangle|^2\leq B_o\|f\|^2.
\end{align*}
The constants $A_o$ and $B_o$ are called the \emph{lower} and \emph{upper frame bounds} of $\Phi$, respectively. For a frame $\Phi=\{\phi_i\}_{i\in[n]}$, the \emph{analysis operator} $T_{\Phi}:\mathbb{C}^k\rightarrow\mathbb{C}^n$ is defined as $T_{\Phi}(f)=\{\langle f,\phi_i\rangle\}_{i\in[n]}$. The analysis operator $T_{\Phi}$ is linear and bounded. The adjoint of $T_{\Phi}$, $T^*_{\Phi}:\mathbb{C}^n\rightarrow\mathbb{C}^k$, is called a \emph{synthesis operator} and defined by $T_{\Phi}(\{\alpha_i\}_{i\in[n]})=\sum\limits_{i\in[n]}\alpha_i \phi_i$. The canonical matrix representation of the synthesis operator is a $k\times n$ matrix such that its $i^{\text{th}}$ column is the vector $\phi_i$ that is  $[T^*_{\Phi}]=[\phi_1\,\,\phi_2\,\,\cdots\,\,\phi_n]$. The \emph{frame operator} $S_{\Phi}:\mathbb{C}^k\rightarrow\mathbb{C}^k$ is the composition of the synthesis operator and analysis operator such that $S_{\Phi}(f)=\sum\limits_{i\in[n]}\langle f,\phi_i\rangle\phi_i$. The operator $S_{\Phi}$ is bounded, linear, self-adjoint, positive, and invertible. The composition $T_{\Phi}T_{\Phi}^*:\mathbb{C}^n\rightarrow\mathbb{C}^n$ is called the \emph{Gramian operator} of $\Phi$. The canonical matrix representation of the Gramian operator is said to be the \emph{Gramian matrix}, defined by $\mathcal{G}_{\Phi}=[\langle \phi_j,\phi_i\rangle]_{i,j\in[n]}$ where $\langle \phi_j,\phi_i\rangle$ is the $(i,j)$-th entry of $\mathcal{G}_{\Phi}$. A frame $\Psi=\{\psi_i\}_{i\in[n]}$ of $\mathbb{C}^k$ is called a \emph{dual frame} of $\Phi$ if for any $f\in\mathbb{C}^k$, $f=\sum\limits_{i\in[n]}\langle f,\phi_i\rangle\psi_i=\sum\limits_{i\in[n]}\langle f,\psi_i\rangle\phi_i$. The sequence $\{S_{\Phi}^{-1}\phi_i\}_{i\in[n]}$ is a also a dual frame of $\Phi$ and it is called the \emph{canonical dual frame} of $\Phi$. We denote the set of all dual frames of a frame $\Phi$ by $D_{\Phi}$. The dual frames which are not canonical dual are called \emph{alternate dual frames}. If $\Psi=\{\psi_i\}_{i\in[n]}\in D_{\Phi}$ then there exist a sequence $\{\mu_i\}_{i\in[n]}\subset\mathbb{C}^k$ such that $\psi_i=S_{\Phi}^{-1}\phi_i+\mu_i$ for $i\in[n]$ and $\sum\limits_{i\in[n]}\langle f,\mu_i\rangle \phi_i=0$ for all $f\in\mathbb{C}^k$. We say two frames $\Phi_1=\{\phi_{1,i}\}_{i\in[n]}$ and $\Phi_2=\{\phi_{2,i}\}_{i\in[n]}$ are \emph{unitary equivalent} if there exists a unitary operator $V:\mathbb{C}^k\rightarrow\mathbb{C}^k$ such that $V(\phi_{1,i})=\psi_{1,i}$ for all $i\in[n]$. Refer to \cite{BL, C, CK1, C2} for the general theory of frames and its applications.

Suppose $\Phi=\{\phi_i\}_{i\in[n]}$ is a frame for $\mathbb{C}^k$ and $\Psi=\{\psi_i\}_{i\in[n]}$ is a dual frame of $\Phi$. Then for $\Lambda\subset[n]$, the error operator $E_{\Phi,\Psi,\Lambda}:\mathbb{C}^k\rightarrow\mathbb{C}^k$ is defined by 
	\begin{align*}
		E_{\Phi,\Psi,\Lambda}(f)=\sum_{i\in\Lambda}\langle f,\phi_i\rangle \psi_i.
	\end{align*}
For $1\leq r<n$, let $\rho^{(r)}_{\Phi,\Psi}:=\max\{\rho(E_{\Phi,\Psi,\Lambda}):\Lambda\subset[n] \text{ and }|\Lambda|=r\}$, where $\rho(E_{\Phi,\Psi,\Lambda})$ is the spectral radius of $E_{\Phi,\Psi,\Lambda}$. A dual frame $\Psi$ of the frame $\Phi$ is called a \emph{spectrally optimal dual frame} (or simply \emph{$SOD$-frame})  for $1$-erasure if $\rho^{(1)}_{\Phi}:=\min\left\{\rho^{(1)}_{\Phi,\widetilde{\Psi}}:\widetilde{\Psi}\in D_{\Phi}\right\}=\rho^{(1)}_{\Phi,\Psi}$. Spectrally optimal dual frames for $r$-erasures are defined inductively. For $1<r<n$, a dual frame $\Psi=\{\psi_i\}_{i\in[n]}$ of a frame $\Phi=\{\phi_i\}_{i\in[n]}$ is called \emph{spectrally optimal dual frame for $r$-erasures} if $\Psi$ is an spectrally optimal dual frame for $(r-1)$-erasures and $\rho^{(r)}_{\Phi}:=\min\{\rho^{(r)}_{\Phi,\widetilde{\Psi}}:\widetilde{\Psi}\in D_{\Phi}\}=\rho^{(r)}_{\Phi,\Psi}$. We denote the set of all $SOD$-frame of $\Phi$ for $r$-erasures by $SOD_{\Phi}^r$.
	  
\begin{remark}\cite{PHM}
    If $\Phi=\{\phi_i\}_{i\in[n]}$ is a frame with dual frame $\Psi=\{\psi_i\}_{i\in[n]}$, then $\rho^{(1)}_{\Phi,\Psi}=\max\{|\langle \psi_i,\phi_i\rangle|:i\in[n]\}$.
\end{remark} 
For the detailed study of optimal dual frames and spectrally optimal dual frames, readers may refer to \cite{NAS, DS1, DS3, LH, LH1, PHM}.

Graph theory is one of the classical branches of mathematics. Along with the deep theoretical aspect of graph theory, it has many real-life applications, which include measuring the connectivity of communication networks, hierarchical clustering, ranking of hyperlinks in search engines, image segmentation, and so on, see \cite{NA, SGT2, SGT3, SGT1}. Let $ G=(V(G),E(G)) $ be a simple graph with the vertex set $ V(G)=\{ v_1, v_2, \dots, v_n\}$  and edge set $ E(G) $. The vertices $ v_i $ and $ v_j $ are said to be \emph{adjacent} if they are connected by an edge and are denoted by  $ v_i \sim v_j $. The number of vertices adjacent to $ v_i $ is known as the \emph{degree of the vertex} $ v_i $ and is denoted by $ d(v_i)$ or simply by $ d_i $. The \emph{degree matrix} of a graph $G$ on $n$ vertices, denoted by $D(G):=diag(d_1,d_2,\ldots,d_n)$. The \emph{adjacency matrix} of $ G $ on $ n $ vertices is $ A(G):=(a_{i,j})_{n\times n} $, where
	\begin{align*}
	a_{i,j}=\begin{cases}
		1, \ \text{for }  v_i\sim v_j \text{ and } i\neq j \\
		0, \text{ elsewhere.}
	\end{cases}
\end{align*}
The matrix $ A(G) $ is simply written as $ A $, if the graph $G$ is clear from the context. The \emph{Laplacian matrix} of a graph $G$ is defined as $L(G):=D(G)-A(G)$. We simply write $ L(G) $ as $ L $ if $ G $ is clear from the context.  The Laplacian matrix is a positive semi-definite matrix. For a simple graph $G$ on $n$ vertices with $k$ components, the rank of $L(G)$ is $n-k$. For a detailed study of matrices associated with simple graphs and other classes of graphs, readers may see \cite{Bapat, SGT6, SGT7, SGT8, SGT9}, and references therein. The transpose of a matrix $ M $ is denoted by $ M^t $. If $ M_1 $ and $ M_2 $ are matrices, then $ M_1 \oplus M_2 $ denotes the block matrix $\left[\begin{array}{cc}
		M_1 & \boldsymbol{0}  \\
		\boldsymbol{0} & M_2
	\end{array}\right]$.

In \cite{D}, authors used Laplacian matrices of graphs to generate finite frames. First, let us see the definition of $L_G(n,k)$-frames.
	\begin{definition}
		Let $G$ be a simple graph on  $n$ vertices with $n-k$ components. Suppose $L$ is the Laplacian matrix of $G$ such that $L=M\, diag(\lambda_1,\lambda_2,\ldots,\lambda_{k},0,\ldots,0)\,M^*$. If $\{e_i\}_{i\in[n]}$ is the standard canonical orthonormal basis of $\mathbb{C}^n$ and $B=diag\left(\sqrt{\lambda_1},\ldots,\sqrt{\lambda_k}\right)M_1^*$, where $M_1$ is a submatrix of $M$ formed by the first $k$ columns then $\{B(e_i)\}_{i\in[n]}$ is a frame for $\mathbb{C}^k$ and it is called an $L_G(n,k)$-frame for $\mathbb{C}^k$.
	\end{definition}
	Next, we have the definition of frames generated by graphs. 
	\begin{definition}
		Let $\Phi=\{\phi_i\}_{i\in [n]}$ be a frame for $\mathbb{C}^k$ with Gramian matrix $\mathcal{G}_{\Phi}$. If $G$ is a graph with Laplacian matrix $L$ such that $\mathcal{G}_{\Phi}=L$, then $\Phi$ is called a frame generated by the graph $G$ for $\mathbb{C}^k$. In short, we call $\Phi$ as a $G(n,k)$-frame for $\mathbb{C}^k$.
	\end{definition}
	In \cite{D}, the authors have shown that every $L_G(n,k)$-frame is a $G(n,k)$-frame. Let us see some results related to frames generated by graphs.
	
	\begin{thm}[\cite{D}]\label{thm2.3}
		If $G$ is a graph, then frames generated by $G$ are unitary equivalent.
	\end{thm}
	The next theorem shows that the frame operator of any $L_G(n,k)$-frame is a diagonal matrix. 
	\begin{thm}[\cite{D}]\label{thm2.4}
		If $G$ is a graph and $\Phi$ is an $L_G(n,k)$-frame then the frame operator $S_{\Phi}$ of the frame $\Phi$ is a diagonal matrix. Further, if $\lambda_1,\lambda_2,\ldots,\lambda_k$ are the non-zero eigenvalues of the Laplacian matrix of $G$ then $S_{\Phi}=diag(\lambda_1,\lambda_2,\ldots,\lambda_k)$.
	\end{thm}
	The following theorem gives the family of dual frames of the frames generated by graphs.
	\begin{thm}\label{thm2.5}\cite{D}
		Let $G$ be a simple graph with components $G_1,G_2,\ldots, G_{m}$ having vertex sets $\{v_1,v_2,\ldots,v_{n_1}\}, \{v_{n_1+1}, v_{n_1+2}, \ldots, v_{n_2}\}, \ldots,\{v_{n_{m-1}+1}, v_{n_{m-1}+2}, \ldots, v_{n_m}(=v_n)\}$, respectively. If $\Phi=\{\phi_i\}_{i\in[n]}$ is a $G(n,n-m)$-frame with the frame operator $S_{\Phi}$, then any dual frame of $\Phi$ is of the form $\{S_{\Phi}^{-1}\phi_i+\nu_1\}_{i=1}^{n_1}\bigcup\{S_{\Phi}^{-1}\phi_i+\nu_2\}_{i=n_1+1}^{n_2}\bigcup\cdots\bigcup\{S_{\Phi}^{-1}\phi_i+\nu_m\}_{i=n_{m-1}+1}^{n}$ where $\nu_1,\nu_2,\ldots,\nu_m$ are arbitrary vectors in $\mathbb{C}^{n-m}$.
	\end{thm}

    If $\{\phi_i\}_{i\in[n]}$ is a frame for $\mathbb{C}^k$ then $\{\phi_i\}_{i\in[n]}$ is called a full spark frame if every subset $\{\phi_i\}_{i\in[n]}$ of cardinality at most k is linearly independent. 
    \begin{thm}\label{thm2.6}\cite{DS2}
		If $G$ is a connected graph and $\Phi=\{\phi_i\}_{i\in[n]}$ is an $L_G(n,n-1)$-frame then $\Phi$ is a full spark frame.
	\end{thm}
In \cite{DS3}, authors studied spectrally optimal dual frames of frames generated by graphs for $1$-erasure and $2$-erasures. Following theorem shows that the canonical dual frames of frames generated by connected graphs are the unique spectrally optimal dual frames for $1$-erasure and $2$-erasures.
\begin{thm}\cite{DS3}\label{thm2.7}
   Suppose $G$ is a connected graph with $n$ vertices. If $\Phi=\{\phi_i\}_{i\in[n]}$ is an $L_G(n,n-1)$-frame for $\mathbb{C}^{n-1}$, then canonical dual frame of $\Phi$ is the unique spectrally optimal dual frame of $\Phi$ for $1$-erasure and $2$-erasures. Further, $\rho^{(2)}_{\Phi}=1$. 
\end{thm}

\begin{thm}\cite{DS3}\label{thm2.8}
    Suppose $G$ is a graph with $n$ vertices and $k>1$ connected components. For any $L_G(n,n-k)$-frame $\Phi$, the canonical dual of $\Psi$ is a non-unique spectrally optimal dual frame of $\Phi$ for $1$-erasure and $2$-erasures. 
\end{thm}

   \textbf{ Note 1:} Suppose $\Phi=\{\phi_i\}_{i\in[n]}$ is an $L_G(n,n-1)$-frame for $\mathbb{C}^{n-1}$, where $G$ is a connected graph on $n$ vertices. Then there exists a diagonal matrix $D=diag(\lambda_1,\lambda_2,\ldots,\lambda_{n-1},0)$ and an orthogonal matrix $M=[M_1, M_2]$, where $M_2=\left[\begin{array}{c}
   \frac{1}{\sqrt{n}}\\
   \frac{1}{\sqrt{n}}\\
   \vdots\\
   \frac{1}{\sqrt{n}}
   \end{array}\right]$ such that the Laplacian matrix $L$ of $G$ is $L=MDM^*$ and $\phi_i=D_1M_1^*(e_i)$ for $i\in[n]$, where $D_1=diag(\sqrt{\lambda_1},\sqrt{\lambda_2},\ldots,\sqrt{\lambda_{n-1}})$ and $\{e_i\}_{i=1}^n$ is an orthonormal basis of $\mathbb{C}^n$. Then by Theorem 2.4, we have
   \begin{align*}
       \|S_{\Phi}^{-1/2}\phi_i\|^2&=\|D_1^{-1}D_1M_1^*(e_i)\|^2\\
       &=\|M_1^*(e_i)\|^2\\
       &=\text{norm of the } i^{\text{th}} \text{row of } M_1^*\\
       &=\text{norm of the } i^{\text{th}} \text{row of } M-\frac{1}{n}\\
       &=1-\frac{1}{n}
   \end{align*}
 For $i\neq j$, we have
 \begin{align*}
     \langle S_{\Phi}^{-1}\phi_i,\phi_j\rangle &=  \langle S_{\Phi}^{-1/2}\phi_i,S_{\Phi}^{-1/2}\phi_j\rangle\\
     &=\langle M_1^*(e_i),M_1^*(e_j)\rangle\\
     &=\text{dot product of }i^{\text{th}} \text{ row of } M_1^* \text{ and } j^{\text{th}} \text{ row of } M_1^*\\
     &=\text{dot product of }i^{\text{th}} \text{ row of } M \text{ and } j^{\text{th}} \text{ row of } M-\frac{1}{n}\\
     &=0-\frac{1}{n}\\
     &=-\frac{1}{n}
 \end{align*}
 \endproof

	\section{Main results}
    We start this section by proving that the spectral radius of the error operator of unitary equivalent frames is equal.
\begin{thm}\label{thm3.1}
    Suppose $\Phi=\{\phi_i\}_{i\in[n]}$ and $\Psi=\{\psi_i\}_{i\in[n]}$ are unitary equivalent frames for $\mathbb{C}^{m}$ with frame operators $S_{\Phi}$ and $S_{\Psi}$, respectively. Then $\rho^{(r)}_{\Phi,S^{-1}_{\Phi}\Phi}=\rho^{(r)}_{\Psi,S^{-1}_{\Psi}\Psi}$ for any $r\in[n-1]$. 
\end{thm}
\proof
Since $\Phi$ and $\Psi$ are unitary equivalent frames, there exists a unitary operator $U:\mathbb{C}^{n}\rightarrow\mathbb{C}^{n}$ such that $\phi_i=U\psi_i$ for $i\in[n]$. For any $f\in\mathbb{C}^n$, we have
\begin{align*}
   S_{\Phi}(f)&=\sum_{i=1}^n\langle f, \phi_i\rangle \phi_i\\
   &=\sum_{i=1}^n\langle f, U\psi_i\rangle U\psi_i\\
    &=U\left(\sum_{i=1}^n\langle U^*f, \psi_i\rangle \psi_i\right)\\
     &=US_{\Psi}U^*(f).
\end{align*}
Thus, $S_{\Phi}=US_{\Psi}U^*$. Let $r\in[n-1]$ be arbitrary. For any $\Lambda\subset[n]$ such that $|\Lambda|=r$, we have
\begin{align*}
   E_{\Phi,S^{-1}_{\Phi}\Phi,\Lambda}(f)&=\sum\limits_{i\in\Lambda}\langle f,\phi_i\rangle S_{\Phi}^{-1}\phi_i\\
   &=\sum\limits_{i\in\Lambda}\langle f,U\psi_i\rangle US_{\Psi}^{-1}U^*U\psi_i\\
     &=U\sum\limits_{i\in\Lambda}\langle U^*f,\psi_i\rangle S_{\Psi}^{-1}\psi_i\\
     &=UE_{\Psi,S^{-1}_{\Psi}\Psi,\Lambda}U^*(f)
\end{align*}
Thus, $E_{\Phi,S^{-1}_{\Phi}\Phi,\Lambda}=UE_{\Psi,S^{-1}_{\Psi}\Psi,\Lambda}U^*$. Hence, we have
\begin{align*}
    \rho^{(r)}_{\Phi,S^{-1}_{\Phi}\Phi}&=\max\{\rho(E_{\Phi,S^{-1}_{\Phi}\Phi,\Lambda}):|\Lambda|=r, \Lambda\subset[n]\}\\
    &=\max\{\rho(UE_{\Psi,S^{-1}_{\Psi}\Psi,\Lambda}U^*):|\Lambda|=r, \Lambda\subset[n]\}\\
     &=\max\{\rho(E_{\Psi,S^{-1}_{\Psi}\Psi,\Lambda}):|\Lambda|=r, \Lambda\subset[n]\}\\
     &=\rho^{(r)}_{\Psi,S^{-1}_{\Psi}\Psi}.
\end{align*}
\endproof

In the next theorem, we study the spectral radius of the error operator for $r$-erasures of frames generated by connected graph.

 \begin{thm}\label{thm4.2}
	If $G$ is a connected graph on $n$ vertices and $\Phi=\{\phi_{i}\}_{i\in[n]}$ is an $L_G(n,n-1)$-frame for $\mathbb{C}^{n-1}$, then $\rho^{(r)}_{\Phi,S^{-1}_{\Phi}\Phi}=1$ for any $r\in[n-1]\setminus\{1\}$. 
\end{thm}	
\proof
By Theorem \ref{thm2.7}, the result is true for $r=2$. Suppose $r\geq3$. Let $L$ be the Laplacian matrix of $G$ with eigenvalues $\lambda_1,\lambda_2,\ldots,\lambda_{n-1},0$. Then there exists an orthogonal matrix $M$ consisting of eigenvectors of $L$ and diagonal matrix $D=diag(\lambda_1,\lambda_2,\ldots ,\lambda_{n-1},0)$ such that $L=MDM^*$ and  $\phi_i=D_1M_1^*(e_i)$ where $D_1=diag(\sqrt{\lambda_1},\sqrt{\lambda_2},\ldots ,\sqrt{\lambda_{n-1}})$, $M_1$ is the matrix formed by columns $1, 2,\ldots,n-1$ of matrix $M$ and $\{e_i\}_{i\in[n]}$ is the standard canonical orthonormal basis of $\mathbb{C}^n$. Note that $[1\,1\,\cdots\,1]^t$ is an eigenvector of $L$ corresponding to the eigenvalue $0$. Hence, the last column of the matrix $M$ is $\left[\frac{1}{\sqrt{n}}\, \frac{1}{\sqrt{n}}\,\cdots\,\frac{1}{\sqrt{n}}\right]^t$. Suppose $S_{\Phi}$ is the frame operator of $\Phi$. Then, by Theorem \ref{thm2.5}, $S_{\Phi}=diag(\lambda_1,\lambda_2,\ldots,\lambda_{n-1})$.

Let $\Lambda\subset[n]$ be arbitrary such that $|\Lambda|=r$. Assume that $\Lambda=\{\alpha_1,\alpha_2,\ldots,\alpha_r\}$. Then, $E_{\Phi,S^{-1}_{\Phi}\Phi,\Lambda}(f)=\sum\limits_{i=1}^r\langle f,\phi_{\alpha_i}\rangle S^{-1}_{\Phi}\phi_{\alpha_i}$. We proceed in the following cases.

\noindent \textbf{Case 1:} Let $r=n$.
 We have $\phi_{\alpha_1}+\phi_{\alpha_2}+\cdots+\phi_{\alpha_{r}}=0$, see the proof of Theorem 4.1 in \cite{DS2}. Thus, $\phi_{\alpha_1}=-(\phi_{\alpha_2}+\cdots+\phi_{\alpha_r})$. Next, we claim that the set $\{\phi_{\alpha_1}-\phi_{\alpha_2},\phi_{\alpha_1}-\phi_{\alpha_3},\ldots,\phi_{\alpha_1}-\phi_{\alpha_r}\}$ is linearly independent. Suppose $c_1(\phi_{\alpha_1}-\phi_{\alpha_2})+c_2(\phi_{\alpha_1}-\phi_{\alpha_3})+\cdots+c_{r-1}(\phi_{\alpha_1}-\phi_{\alpha_r})=0$. Then $\phi_{\alpha_1}=\frac{c_1}{c_1+\cdots+c_{r-1}}\phi_{\alpha_2}+\cdots+\frac{c_{r-1}}{c_1+\cdots+c_{r-1}}\phi_{\alpha_r}$. Thus, we have 
\begin{align*}
   \left( \frac{c_1}{c_1+\cdots+c_{r-1}}+1\right)\phi_{\alpha_2}+\cdots+\left(\frac{c_{r-1}}{c_1+\cdots+c_{r-1}}+1\right)\phi_{\alpha_r}=0
\end{align*}
By Theorem \ref{thm2.6}, $\Phi$ is a full spark frame. Thus, $\{\phi_{\alpha_2},\ldots,\phi_{\alpha_r}\}$ is linearly independent. Hence, we have
\begin{align*}
    2c_1+c_2+\cdots+c_{r-1}&=0\\
    c_1+2c_2+\cdots+c_{r-1}&=0\\
    \vdots\quad\quad\quad&\\
    c_1+c_2+\cdots+2c_{r-1}&=0
\end{align*}
This gives $c_1=c_2=\cdots=c_{r-1}=0$. Therefore, the set $\{\phi_{\alpha_1}-\phi_{\alpha_2},\phi_{\alpha_1}-\phi_{\alpha_3},\ldots,\phi_{\alpha_1}-\phi_{\alpha_r}\}$ is linearly independent. Since $S_{\Phi}$ is an invertible operator,  $\{S_{\Phi}^{-1}\phi_{\alpha_1}-S_{\Phi}^{-1}\phi_{\alpha_2},S_{\Phi}^{-1}\phi_{\alpha_1}-S_{\Phi}^{-1}\phi_{\alpha_3},\ldots,S_{\Phi}^{-1}\phi_{\alpha_1}-S_{\Phi}^{-1}\phi_{\alpha_r}\}$ is linearly independent.

\noindent \textbf{Case 2:} Let $r<n$. Since $\Phi$ is a full spark frame, $\{\phi_{\alpha_1},\ldots, \phi_{\alpha_r} \}$ is linearly independent. Hence, $\{\phi_{\alpha_1}-\phi_{\alpha_2},\phi_{\alpha_1}-\phi_{\alpha_3},\ldots,\phi_{\alpha_1}-\phi_{\alpha_r}\}$ is linearly independent. Thus, $\{S_{\Phi}^{-1}\phi_{\alpha_1}-S_{\Phi}^{-1}\phi_{\alpha_2},S_{\Phi}^{-1}\phi_{\alpha_1}-S_{\Phi}^{-1}\phi_{\alpha_3},\ldots,S_{\Phi}^{-1}\phi_{\alpha_1}-S_{\Phi}^{-1}\phi_{\alpha_r}\}$ is linearly independent.

From the above two cases, we have $\{S_{\Phi}^{-1}\phi_{\alpha_1}-S_{\Phi}^{-1}\phi_{\alpha_2},S_{\Phi}^{-1}\phi_{\alpha_1}-S_{\Phi}^{-1}\phi_{\alpha_3},\ldots,S_{\Phi}^{-1}\phi_{\alpha_1}-S_{\Phi}^{-1}\phi_{\alpha_r}\}$ is linearly independent.  Then, for any $i\in\{2,3,\ldots,r\}$, we have
\begin{align*}
    E_{\Phi,S^{-1}_{\Phi}\Phi,\Lambda}(S_{\Phi}^{-1}\phi_{\alpha_1}-S_{\Phi}^{-1}\phi_{\alpha_i})&=\sum_{j\in[r]}\langle S_{\Phi}^{-1}\phi_{\alpha_1}-S_{\Phi}^{-1}\phi_{\alpha_i}, \phi_{\alpha_j}\rangle S^{-1}_{\Phi} \phi_{\alpha_j}\\
    &=\sum_{j\in[r]}(\langle S_{\Phi}^{-1}\phi_{\alpha_1}, \phi_{\alpha_j}\rangle-\langle S_{\Phi}^{-1}\phi_{\alpha_i}, \phi_{\alpha_j}\rangle) S^{-1}_{\Phi} \phi_{\alpha_j}\\
    &=\left(\left(1-\frac{1}{n}\right)-\left(-\frac{1}{n}\right)\right)S_{\Phi}^{-1}\phi_{\alpha_1}+\left(\left(-\frac{1}{n}\right)-\left(1-\frac{1}{n}\right)\right)S_{\Phi}^{-1}\phi_{\alpha_i}\\
    &+\sum_{j\in[r],j\neq 1, i}\left(\left(-\frac{1}{n}\right)-\left(-\frac{1}{n}\right)\right)S_{\Phi}^{-1}\phi_{\alpha_j}\\
    &=S^{-1}_{\Phi}\phi_{\alpha_1}-S^{-1}_{\Phi}\phi_{\alpha_i}.
\end{align*}
Thus, $S_{\Phi}^{-1}\phi_{\alpha_1}-S_{\Phi}^{-1}\phi_{\alpha_2},S_{\Phi}^{-1}\phi_{\alpha_1}-S_{\Phi}^{-1}\phi_{\alpha_3},\ldots,S_{\Phi}^{-1}\phi_{\alpha_1}-S_{\Phi}^{-1}\phi_{\alpha_r}$ are the eigenvectors of $E_{\Phi,S^{-1}_{\Phi}\Phi,\Lambda}$ corresponding to the eigenvalue $1$. Hence, $1$ an eigenvalue of $E_{\Phi,S^{-1}\Phi,\Lambda}$ with multiplicity $r-1$. If $r=n$ then rank of the error operator $E_{\Phi,S^{-1}_{\Phi}\Phi,\Lambda}$ is $r-1$. Hence, the operator $E_{\Phi,S^{-1}\Phi,\Lambda}$ has $r-1$ non-zero eigenvalues and thus $1$ is the only non-zero eigenvalue of $E_{\Phi,S^{-1}_{\Phi}\Phi,\Lambda}$. If $r<n$, then the rank of the error operator $E_{\Phi,S^{-1}_{\Phi}\Phi,\Lambda}$ is $r$. Hence, the operator $E_{\Phi,S^{-1}\Phi,\Lambda}$ has $r$ non-zero eigenvalues and $S_{\Phi}^{-1}\phi_{\alpha_1}+\cdots+S_{\Phi}^{-1}\phi_{\alpha_r}\neq 0$ as the set $\{S_{\Phi}^{-1}\phi_{\alpha_1},\ldots,S_{\Phi}^{-1}\phi_{\alpha_r}\}$ is linearly independent. We have
\begin{align*}
    E_{\Phi,S_{\Phi}^{-1}\Phi,\Lambda}(S_{\Phi}^{-1}\phi_{\alpha_1}&+\cdots+S_{\Phi}^{-1}\phi_{\alpha_r})=\sum_{i=1}^r\langle S_{\Phi}^{-1}\phi_{\alpha_1}+\cdots+S_{\Phi}^{-1}\phi_{\alpha_r},\phi_{\alpha_i}\rangle S_{\Phi}^{-1}\phi_{\alpha_i}\\
    &=\sum_{i=1}^r(\langle S_{\Phi}^{-1}\phi_{\alpha_1},\phi_{\alpha_i}\rangle + \cdots+ \langle S_{\Phi}^{-1}\phi_{\alpha_i},\phi_{\alpha_i}\rangle+\cdots+\langle S_{\Phi}^{-1}\phi_{\alpha_r},\phi_{\alpha_i}\rangle)S_{\Phi}^{-1}\phi_{\alpha_i}\\
    &=\sum_{i=1}^r\left( \left(-\frac{1}{n}\right)+ \cdots+ \left(1-\frac{1}{n}\right)+\cdots+\left(-\frac{1}{n}\right)\right)S_{\Phi}^{-1}\phi_{\alpha_i}\\
    &=\sum_{i=1}^r\left(1-\frac{r}{n}\right)S_{\Phi}^{-1}\phi_{\alpha_i}\\
     &=\left(1-\frac{r}{n}\right)\sum_{i=1}^rS_{\Phi}^{-1}\phi_{\alpha_i}.
\end{align*}
Thus, $1-\frac{r}{n}$ is a non-zero eigenvalue of $E_{\Phi,S_{\Phi}^{-1}\Phi,\Lambda}$ if $r<n$. Hence, we have
\begin{align*}
    \sigma(E_{\Phi,S_{\Phi}^{-1}\Phi,\Lambda})=\begin{cases}
        \left\{1,1-\frac{r}{n},0\right\} \text{ if } r<n\\
        \{1,0\}  \text{ if } r=n.
    \end{cases}
\end{align*}
Thus, $\rho(E_{\Phi,S_{\Phi}^{-1}\Phi,\Lambda})=1$. Hence, $\rho^{(r)}_{\Phi, S^{-1}_{\Phi}\Phi}=\max\{\rho(E_{\Phi,S_{\Phi}^{-1}\Phi,\Lambda}):\Lambda\subset[n], |\Lambda|=r\}=1$.
\endproof

\begin{thm}\label{thm3.3}
If $G$ is a graph with $n$ vertices and $k>1$ connected components, then there exists an $L_{G}(n,n-k)$-frame $\Phi$ for $\mathbb{C}^{n-k}$ such that $\rho^{(r)}_{\Phi, S_{\Phi}^{-1}\Phi}=1$ for any $r\in[n-1]\setminus\{1\}$. Further, $S_{\Phi}^{-1}\Phi$ is a non-unique SOD-frame of $\Phi$ for $r$-erasures for any $r\in[n-1]$.
\end{thm}
\proof
Let $G_1,G_2,\ldots,G_k$ be the connected components of $G$ with $|V(G_i)|=n_i$ for $i\in[k]$. Let the vertex set of $G_i$ be $\{v_{l_{i-1}+1},v_{l_{i-1}+2},\ldots,v_{l_i}\}$ for $i\in[k]$ where $l_0=0$ and $l_i=n_1+n_2+\cdots+n_i$  for $i\in[k]$. If $L_i$ is the Laplacian matrix of $G_i$, then the Laplacian matrix of $G$ is $L=L_1\oplus L_2\oplus\cdots\oplus L_k$. Suppose spectral decomposition of $L_i=M_iD_iM_i^*$ where $D_i=diag(\lambda_1^i,\ldots,\lambda_{n_i-1}^i,0)$ and $M_i$ is an orthogonal matrix consist of eigenvectors of $L_i$. Then $L=MDM^*$ where $M=M_1\oplus\cdots\oplus M_k$ and $D=D_1\oplus\cdots\oplus D_k$. Let $\widetilde{D_i}=diag\left(\sqrt{\lambda_1^i},\ldots,\sqrt{\lambda_{n_i-1}^i}\right)$ and $\widetilde{M_i}$ is formed by columns $1,2,\ldots,n-1$ of matrix $M_i$ for $i\in[k]$. Then $\widetilde{D_i}\widetilde{M_i}^*$ is the synthesis operator of frame, say $\Phi_i=\{\phi_{i,j}\}_{j\in[n_i]}$ for the space $\mathbb{C}^{n_i-1}$. Let $S_{\Phi_i}$ be the frame operator of $\Phi_i$ for $i\in[k]$, $\widetilde{D}=\widetilde{D_1}\oplus\cdots\oplus\widetilde{D_k}$ and $\widetilde{M}
=\widetilde{M_1}\oplus\cdots\oplus\widetilde{M_k}$. Then $\Phi=\{\phi_i\}_{i\in[n]}=\{\widetilde{D}\widetilde{M}^*(e_i)\}_{i\in[n]}$ is an $L_G(n,n-k)$-frame for $\mathbb{C}^{n-k}$. Also, the synthesis operator $T_{\Phi}^*$ of $\Phi$ is
\begin{align*}
	[T_{\Phi}^*]&=[T_{\Phi_1}^*]\oplus[T_{\Phi_2}^*]\oplus\cdots\oplus[T_{\Phi_k}^*]\\
	&=\left[\begin{array}{cccccccccc}
		\phi_{1,1}  &\cdots & \phi_{1,n_1} & 0 & \cdots & 0 & \cdots & 0 & \cdots & 0\\
		0  &\cdots & 0 & \phi_{2,1} & \cdots & \phi_{2,n_2} & \cdots & 0 & \cdots & 0\\
		\vdots  &\ddots & \vdots & \vdots & \ddots & \vdots & \ddots & \vdots & \ddots & \vdots \\
		0  &\cdots & 0 & 0 & \cdots & 0 & \cdots & \phi_{k,1} & \cdots & \phi_{k,n_k}
	\end{array}\right]\\
	&=\left[\begin{array}{cccc}
		\phi_{1}  & \phi_{2} & \cdots & \phi_{n} 
	\end{array}\right]\\
\end{align*}
Let $S_{\Phi}$ be the frame operator of $\Phi$. Then $[S_{\Phi}]=[S_{\Phi_1}]\oplus\cdots\oplus [S_{\Phi_k}]$. By Theorem \ref{thm2.8}, $S_{\Phi}^{-1}\Phi$ is a non-unique SOD-frame of $\Phi$ for $1$-erasure and $2$-erasures. For any $\Lambda\subset[n]$ such that $|\Lambda|=r\neq1$, $E_{\Phi,S^{-1}_{\Phi}\Phi,\Lambda}(f)=\sum\limits_{i\in\Lambda}\langle f,\phi_i\rangle S_{\Phi}^{-1}\phi_i$. Let $\Lambda=\{\alpha_1,\alpha_2, \ldots, \alpha_r\}$ such that $\{\alpha_1,\ldots,\alpha_{m_1}\}\subset V(G_1)$, $\{\alpha_{m_1+1},\ldots,\alpha_{m_2}\}\subset V(G_2)$, $\ldots$, $\{\alpha_{m_{k-1}+1},\ldots,\alpha_{m_k}\}\subset V(G_k)$ where $m_k=r$. 

Let $m_0=0$ and $\widetilde{\Lambda}=\{i\in[k]:\{\alpha_{m_{i-1}+1},\ldots,\alpha_{m_i}\}\neq\emptyset\}$. Assume $\Lambda_0=\{i\in\widetilde{\Lambda}:|\{\alpha_{m_{i-1}+1}, \ldots,\alpha_{m_i}\}|=1\neq |V(G_i)|\}$, $\Lambda_1=\{i\in\widetilde{\Lambda}:|\{\alpha_{m_{i-1}+1}, \ldots,\alpha_{m_i}\}|=1=|V(G_i)|\}$, $\Lambda_2=\{i\in\widetilde{\Lambda}\setminus\Lambda_1:|\{\alpha_{m_{i-1}+1}, \ldots,\alpha_{m_i}\}|=|V(G_i)|\}$ and $\Lambda_3=\widetilde{\Lambda}\setminus (\Lambda_0\cup\Lambda_1\cup\Lambda_2)$. Since $G_j$ is connected, $\{\phi_{j,1},\ldots,\phi_{j,n_j}\}$ is a full spark frame for $\mathbb{C}^{n_j-1}$ for $j\in[k]$. Then, the rank of the error operator $E_{\Phi,S^{-1}_{\Phi}\Phi,\Lambda}$ is $|\Lambda_0|+\sum\limits_{i\in\Lambda_2}(n_i-1)+\sum\limits_{i\in\Lambda_3}|\{\alpha_{m_{i-1}+1}, \ldots,\alpha_{m_i}\}|$.

Let $j\in\Lambda_2$. Now, $\{\phi_{j,1},\ldots,\phi_{j,n_j}\}$ is a full spark frame for $\mathbb{C}^{n_j-1}$ and hence any subset of $\{\phi_{j,1},\ldots,\phi_{j,n_j}\}$ of cardinality at most $n_j-1$ is linearly independent. Thus, any subset of $\{\phi_{\alpha_{m_{j-1}+1}},\phi_{\alpha_{m_{j-1}+2}},\ldots,\phi_{\alpha_{m_j}}\}$ of cardinality at most $n_j-1$ is linearly independent. Also, we have $\phi_{j,1}+\cdots+\phi_{j,n_j}=0$, see the proof of Theorem 4.1 in \cite{DS2}. Then, $\phi_{\alpha_{m_{j-1}+1}}+\phi_{\alpha_{m_{j-1}+2}}+\cdots+\phi_{\alpha_{m_j}}=0$. Hence, the set $\{\phi_{\alpha_{m_j-1}+1}- \phi_{\alpha_{m_j-1}+2},\ldots,\phi_{\alpha_{m_j-1}+1}-\phi_{\alpha_{m_j}}\}$ is linearly independent. This gives that the set $\{S^{-1}_{\Phi}\phi_{\alpha_{m_j-1}+1}- S^{-1}_{\Phi}\phi_{\alpha_{m_j-1}+2},\ldots,S^{-1}_{\Phi}\phi_{\alpha_{m_j-1}+1}-S^{-1}_{\Phi}\phi_{\alpha_{m_j}}\}$ is linearly independent.

If $j\in\Lambda_3$, then $\{\phi_{\alpha_{m_{j-1}+1}},\phi_{\alpha_{m_{j-1}+2}},\ldots,\phi_{\alpha_{m_j}}\}$ is linearly independent. Thus, the set $\{\phi_{\alpha_{m_j-1}+1}- \phi_{\alpha_{m_j-1}+2},\ldots,\phi_{\alpha_{m_j-1}+1}-\phi_{\alpha_{m_j}}\}$ is linearly independent. Therefore, $\{S^{-1}_{\Phi}\phi_{\alpha_{m_j-1}+1}- S^{-1}_{\Phi}\phi_{\alpha_{m_j-1}+2},\ldots,S^{-1}_{\Phi}\phi_{\alpha_{m_j-1}+1}-S^{-1}_{\Phi}\phi_{\alpha_{m_j}}\}$ is linearly independent.

For any $j\in\Lambda_0$, $\{\phi_{\alpha_{m_j}}\}$ is linearly independent that is $\phi_{\alpha_{m_j}}\neq 0$. Thus, $S^{-1}_{\Phi} \phi_{\alpha_{m_j}}\neq 0$.

For any $l\in(\Lambda_2\cup\Lambda_3)$ and $j\in\{2,3,\ldots,m_l-m_{l-1}\}$, we have 
\begin{align*}
    E_{\Phi,S^{-1}_{\Phi}\Phi,\Lambda}&(S_{\Phi}^{-1}\phi_{\alpha_{m_{l-1}+1}}-S_{\Phi}^{-1}\phi_{\alpha_{m_{l-1}+j}})\\
   &=\sum_{i\in\widetilde{\Lambda}}\sum_{q=1}^{m_i-m_{i-1}}\langle S_{\Phi}^{-1}\phi_{\alpha_{m_{l-1}+1}}-S_{\Phi}^{-1}\phi_{\alpha_{m_{l-1}+j}}, \phi_{\alpha_{m_{i-1}+q}}\rangle S_{\Phi}^{-1}\phi_{\alpha_{m_{i-1}+q}}\\
     &=\sum_{q=1}^{m_l-m_{l-1}}\langle S_{\Phi}^{-1}\phi_{\alpha_{m_{l-1}+1}}-S_{\Phi}^{-1}\phi_{\alpha_{m_{l-1}+j}}, \phi_{\alpha_{m_{l-1}+q}}\rangle S_{\Phi}^{-1}\phi_{\alpha_{m_{l-1}+q}}\\
    &=\left(\left(1-\frac{1}{n_l}\right)-\left(-\frac{1}{n_l}\right)\right)S_{\Phi}^{-1}\phi_{\alpha_{m_{l-1}+1}}+\left(\left(-\frac{1}{n_l}\right)-\left(1-\frac{1}{n_l}\right)\right)S_{\Phi}^{-1}\phi_{\alpha_{m_{l-1}+j}}\\
      &=S_{\Phi}^{-1}\phi_{\alpha_{m_{l-1}+1}}-S_{\Phi}^{-1}\phi_{\alpha_{m_{l-1}+j}}
\end{align*}
For $l\in\Lambda_3$, we have
\begin{align*}
      E_{\Phi,S^{-1}_{\Phi}\Phi,\Lambda}&(S_{\Phi}^{-1}\phi_{\alpha_{m_{l-1}+1}}+\cdots+S_{\Phi}^{-1}\phi_{\alpha_{m_l}})\\
&=\sum_{i\in\widetilde{\Lambda}}\sum_{q=1}^{m_i-m_{i-1}}\langle S_{\Phi}^{-1}\phi_{\alpha_{m_{l-1}+1}}+\cdots+S_{\Phi}^{-1}\phi_{\alpha_{m_l}}, \phi_{\alpha_{m_{i-1}+q}}\rangle S_{\Phi}^{-1}\phi_{\alpha_{m_{i-1}+q}}\\
       &=\sum_{q=1}^{m_l-m_{l-1}}\langle S_{\Phi}^{-1}\phi_{\alpha_{m_{l-1}+1}}+\cdots+S_{\Phi}^{-1}\phi_{\alpha_{m_l}}, \phi_{\alpha_{m_{l-1}+q}}\rangle S_{\Phi}^{-1}\phi_{\alpha_{m_{l-1}+q}}\\
       &=\sum_{q=1}^{m_i-m_{i-1}}\left(1-\frac{m_l}{n_l}\right) S_{\Phi}^{-1}\phi_{\alpha_{m_{i-1}+q}}\\
        &=\left(1-\frac{m_l}{n_l}\right)(S_{\Phi}^{-1}\phi_{\alpha_{m_{l-1}+1}}+\cdots+S_{\Phi}^{-1}\phi_{\alpha_{m_l}})\\
\end{align*}

Further, for any $l\in\Lambda_0$, we have
\begin{align*}
      E_{\Phi,S^{-1}_{\Phi}\Phi,\Lambda}(S_{\Phi}^{-1}\phi_{\alpha_{m_{l-1}+1}})&= E_{\Phi,S^{-1}_{\Phi}\Phi,\Lambda}(S_{\Phi}^{-1}\phi_{\alpha_{m_{l}}})\\
      &=\sum_{i\in\widetilde{\Lambda}}\sum_{q=1}^{m_i-m_{i-1}}\langle S_{\Phi}^{-1}\phi_{\alpha_{m_l}}, \phi_{\alpha_{m_{i-1}+q}}\rangle S_{\Phi}^{-1}\phi_{\alpha_{m_{i-1}+q}}\\
      &=\langle S_{\Phi}^{-1}\phi_{\alpha_{m_l}}, \phi_{\alpha_{m_{l}}}\rangle S_{\Phi}^{-1}\phi_{\alpha_{m_{l}}}\\
      &=\left(1-\frac{1}{n_l}\right)S_{\Phi}^{-1}\phi_{\alpha_{m_l}}.
\end{align*}
Therefore, the eigenvalues of $E_{\Phi,S^{-1}_{\widetilde{\Phi}}\widetilde{\Phi},\Lambda}$ are $1-\frac{m_l}{n_l}$ with multiplicity $1$ for $l\in\Lambda_3$, $1-\frac{1}{n_l}$ with multiplicity $1$ for $l\in\Lambda_0$ and $1$ with multiplicity $m_l-m_{l-1}-1$ for $l\in\Lambda_2\cup\Lambda_3$. Hence, $\rho(E_{\Phi,S^{-1}_{\Phi}\Phi,\Lambda})=
    1, \text{ if } \Lambda_2\cup\Lambda_3\neq\emptyset$ and $\rho(E_{\Phi,S^{-1}_{\Phi}\Phi,\Lambda})\leq
    1, \text{ otherwise }$. Therefore, $\rho^{(r)}_{\Phi,S^{-1}_{\Phi}\Phi}=\max\{\rho(E_{\Phi,S^{-1}_{\Phi}\Phi,\Lambda}):|\Lambda|=r, \Lambda\subset[n]\}=1$ for any $r\in[n-1]\setminus\{1\}$. 

Let $\Psi=\{\psi_i\}_{i\in[n]}$ be any dual frame of $\Phi$. If $G$ is a null graph, then there is nothing to prove. Hence, assume that $G$ is not a null graph. Thus, there exists $j_o\in[k]$ such that $|V(G_{j_o})|\geq 2$. Let $p,q\in V(G_{j_o}) $. Choose $\Lambda\subset[n]$ such that $|\Lambda|=r$ and $p,q\in\Lambda$. By Theorem \ref{thm2.5}, there exists $\nu\in\mathbb{C}^{n-k}$ such that $\psi_i=S_{\Phi}^{-1}\phi_i+\nu$ for $i\in V(G_{j_o})$. Now $E_{\Phi,\Psi,\Lambda}(f)=\sum\limits_{i\in\Lambda}\langle f,\phi_i\rangle \psi_i$. Then by using Note 1, we have
\begin{align*}
   E_{\Phi,\Psi,\Lambda}(\psi_p-\psi_q)&=\sum_{i\in\Lambda} \langle \psi_p-\psi_q,\phi_i\rangle \psi_i\\
   &=\sum_{i\in\Lambda} \langle S_{\Phi}^{-1}\phi_p-S_{\Phi}^{-1}\phi_q,\phi_i\rangle \psi_i\\
    &=\langle S_{\Phi}^{-1}\phi_p-S_{\Phi}^{-1}\phi_q,\phi_p\rangle \psi_p+\langle S_{\Phi}^{-1}\phi_p-S_{\Phi}^{-1}\phi_q,\phi_q\rangle \psi_q\\
    &=\left(1-\frac{1}{n_{j_o}}+\frac{1}{n_{j_o}}\right)\psi_p+\left(-\frac{1}{n_{j_o}}-1+\frac{1}{n_{j_o}}\right)\psi_q\\
    &=\psi_p-\psi_q.
\end{align*}
Thus, $1\in\sigma(E_{\Phi,\Psi,\Lambda})$ that is $\rho(E_{\Phi,\Psi,\Lambda})\geq 1$. Hence, $\rho_{\Phi,\Psi}^{(r)}\geq 1=\rho_{\Phi,S^{-1}_{\Phi}\Phi}^{(r)}$. By Theorem \ref{thm2.6}, $S_{\Phi}^{-1}\Phi$ is an $SOD$-frame of $\Phi$ for $1$-erasure and $\rho_{\Phi,S^{-1}_{\Phi}\Phi}^{(r)}\leq \rho_{\Phi,\Psi}^{(r)}$ for any dual frame $\Psi$ of $\Phi$ and for any $r\in[n]$, therefore, $S_{\Phi}^{-1}\Phi$ is an $SOD$-frame for $\Phi$ for $r$-erasures.

Let $V(G_1)=\{1,\ldots,\alpha\}$. Choose a dual frame $\widetilde{\Psi}=\{\widetilde{\psi_i}\}_{i\in[n]}=\{S_{\Phi}^{-1}\phi_i+\nu\}_{i\in[\alpha]}\cup\{S_{\Phi}^{-1}\phi_i\}_{i\in[n\setminus[\alpha]]}$ where $\nu=\left[\begin{array}{c}
     v_1\\
     v_2\\
     \vdots\\
     v_{n-k}
\end{array}\right]$ such that $v_i=\begin{cases}
    1, \text{ if } i=\alpha\\
    0, \text{ otherwise}
\end{cases}$. Then, we have
\begin{align*}
    \rho^{(1)}_{\Phi,\widetilde{\Psi}}&=\max\{|\langle\widetilde{\psi_i}, \phi_i\rangle|:i\in[n]\}\\
    &=\max(\{|\langle S_{\Phi}^{-1}\phi_i+\nu, \phi_i\rangle|:i\in[\alpha]\}\cup\{|\langle S_{\Phi}^{-1}\phi_i, \phi_i\rangle|:i\in[n]\setminus[\alpha]\})\\
     &=\max(\{|\langle S_{\Phi}^{-1}\phi_i, \phi_i\rangle|:i\in[\alpha]\}\cup\{|\langle S_{\Phi}^{-1}\phi_i, \phi_i\rangle|:i\in[n]\setminus[\alpha]\})\\
&=\rho^{(1)}_{\Phi,S_{\Phi}^{-1}\Phi}
\end{align*}
Therefore, $\widetilde{\Psi}$ is an $SOD$-frame of $\Phi$ for $1$-erasure. Let $\Lambda\subset[n]$ be arbitrary such that $|\Lambda|=r\geq 2$. Let $\Lambda=\Lambda_1\cup\Lambda_2\cup\cdots\cup\Lambda_k$ such that $\Lambda_i=\{\alpha_{l_{i-1}+1},\alpha_{l_{i-1}+2},\ldots,\alpha_{l_{i}}\}\subset V(G_i)$ for $i\in[k]$ where $l_0=0$ and $l_k=r$. 

If $|\{\alpha_1,\alpha_2,\ldots,\alpha_{l_1}\}|=\emptyset$, then, $E_{\Phi,\widetilde{\Psi},\Lambda}=E_{\Phi,S_{\Phi}^{-1}\Phi,\Lambda}$. Thus, $\rho^{(r)}_{\Phi,\widetilde{\Psi}}=\rho^{(r)}_{\Phi,S_{\Phi}^{-1}\Phi}=1$.

If $|\{\alpha_1,\alpha_2,\ldots,\alpha_{l_1}\}|\neq\emptyset$, then $E_{\Phi,\widetilde{\Psi},\Lambda}(f)=\sum\limits_{i=1}^{l_1}\langle f,\phi_i\rangle (S_{\Phi}^{-1}\phi_i+\nu)+\sum\limits_{j=2}^k\,\sum\limits_{i=l_{j-1}+1}^{l_j}\langle f, \phi_i\rangle S^{-1}_{\Phi}\phi_i$ and $E_{\Phi,S_{\Phi}^{-1}\Phi,\Lambda}(f)=\sum\limits_{j=1}^k\,\sum\limits_{i=l_{j-1}+1}^{l_j}\langle f, \phi_i\rangle S^{-1}_{\Phi}\phi_i$. Then, the matrix representation of $E_{\Phi,S_{\Phi}^{-1}\Phi,\Lambda}$ is of the form $\left[\begin{array}{cc}
  E_1   &  0\\
  0   &  E_2
\end{array}\right]$ where $E_1$ is a matrix of order $\alpha-1\times \alpha-1$ and $E_2$ is a matrix of order $(n-k)-(\alpha-1)\times(n-k)-(\alpha-1)$. Thus, matrix representation of $E_{\Phi, \widetilde{\Psi},\Lambda}$  is $\left[\begin{array}{cc}
  E_1   &  0\\
  E_3   &  E_2
\end{array}\right]$. Therefore, characterstic polynomial of $E_{\Phi, \widetilde{\Psi},\Lambda}$ is same as the characterstic polynomial of $E_{\Phi,S_{\Phi}^{-1}\Phi,\Lambda}$. Hence, $\sigma(E_{\Phi,S_{\Phi}^{-1}\Phi,\Lambda})=\sigma(E_{\Phi, \widetilde{\Psi},\Lambda})$. Thus, we have $\rho^{(r)}_{\Phi,S^{-1}_{\Phi}\Phi}=\rho^{(r)}_{\Phi,\widetilde{\Psi}}=1$. Therefore, $\widetilde{\Psi}$ is a non-canonical spectrally optimal dual frame of $\Phi$ for $r$-erasures. Thus, $S_{\Phi}^{-1}\Phi$ is the non-unique $SOD$-frame of $\Phi$ for $r$-erasures.  
\endproof

\begin{thm}\label{thm4.3}
	Suppose $G$ is a graph with $n$ vertices and $k>1$ connected components. If $\Phi=\{\phi_i\}_{i\in[n]}$ is an $L_G(n,n-k)$-frame for $\mathbb{C}^{n-k}$, then $\rho^{(r)}_{\Phi, S_{\Phi}^{-1}\Phi}=1$ for any $r\in[n-1]\setminus\{1\}$. 
\end{thm}	
\proof
By Theorem \ref{thm3.3}, there exists an $L_G(n,n-k)$-frame $\Psi=\{\psi_i\}_{i\in[n]}$ such that $\rho^{(r)}_{\Psi, S_{\Psi}^{-1}\Psi}=1$ for any $r\in[n-1]\setminus\{1\}$. By Theorem \ref{thm2.3}, the frames $\Phi$ and $\Psi$ are unitary equivalent. Using Theorem \ref{thm3.1}, $\rho^{(r)}_{\Phi,S^{-1}_{\Phi}\Phi}=\rho^{(r)}_{\Psi,S^{-1}_{\Psi}\Psi}=1$ for any $r\in[n-1]\setminus\{1\}$.
\endproof
Next theorem shows that the canonical dual frames of frames generated by connected graphs are the unique $SOD$-frame for $r$-erasures.
 \begin{thm}
     If $G$ is a connected graph on $n$-vertices and $\Phi=\{\phi_i\}_{i\in[n]}$ is an $L_G(n,n-1)$-frame for $\mathbb{C}^{n-1}$ then for any $r\in[n-1],$ $S_{\Phi}^{-1}\Phi$ is the unique $SOD$-frame of $\Phi$ for $r$-erasures.
 \end{thm}
 \proof
By Theorem \ref{thm2.7}, the result is true for $r=1,2$. Thus, assume $r\geq 3$. Let $\Psi=\{\psi_i\}_{i\in[n]}$ be any dual frame of $\Phi$. By Theorem \ref{thm2.5}, there exists a vector $\nu\in\mathbb{C}^{n-1}$ such that $\psi_i=S_{\Phi}^{-1}\phi_i+\nu$ for $i\in[n]$. 

Let $\Lambda\subset[n]$ be arbitrary such that $|\Lambda|=r$. Suppose $i,j\in\Lambda$. Then, we have
\begin{align*}
    E_{\Phi,\Psi,\Lambda}(\psi_i-\psi_j)&=\sum_{q\in\Lambda}\langle\psi_i-\psi_j, \phi_q\rangle \psi_q\\
&=\sum_{q\in\Lambda}\langle S_{\Phi}^{-1}\phi_i-S_{\Phi}^{-1}\phi_j, \phi_q\rangle \psi_q\\
&=\sum_{{q\in\Lambda},{q\neq i,j}}(\langle S_{\Phi}^{-1}\phi_i, \phi_q\rangle-\langle S_{\Phi}^{-1}\phi_j, \phi_q\rangle) \psi_q+(\langle S_{\Phi}^{-1}\phi_i, \phi_i\rangle-\langle S_{\Phi}^{-1}\phi_j, \phi_i\rangle) \psi_i\\
&\quad\quad\quad\quad+ (\langle S_{\Phi}^{-1}\phi_i, \phi_j\rangle-\langle S_{\Phi}^{-1}\phi_j, \phi_j\rangle) \psi_j\\
&=\sum_{{q\in\Lambda},{q\neq i,j}}\left(\left(-\frac{1}{n}\right)-\left(-\frac{1}{n}\right)\right) \psi_q+\left(\left(1-\frac{1}{n}\right)-\left(-\frac{1}{n}\right)\right) \psi_i\\
&\quad\quad\quad\quad+ \left(\left(-\frac{1}{n}\right)-\left(1-\frac{1}{n}\right)\right)\psi_j\\
&=\psi_i-\psi_j.
\end{align*}
Therefore, $1\in\sigma( E_{\Phi,\Psi,\Lambda})$. Thus, $\rho^{(r)}_{\Phi,\Psi}\geq 1=\rho^{(r)}_{\Phi,S^{-1}_{\Phi}\Phi}$. Since, $S_{\Phi}^{-1}\Phi$ is the unique $SOD$-frame for $1$-erasure and $\rho^{(r)}_{\Phi,S^{-1}_{\Phi}\Phi}\leq\rho^{(r)}_{\Phi,\Psi}$ for any dual frame $\Psi$ of $\Phi$, thus, $S_{\Phi}^{-1}\Phi$ is the unique $SOD$-frame for $r$-erasures.
 \endproof

In \cite{DS3}, it is shown that canonical dual of frames generated by disconnected graph is a non-unique $SOD$-frame for $1$-erasure and $2$-erasures. In the following theorem, we show that the result is true for $r$-erasures for any $r$.

 \begin{thm}
     If $G$ is a graph on $n$-vertices with $k>1$ connected components and $\Phi=\{\phi_i\}_{i\in[n]}$ is an $L_G(n,n-k)$-frame for $\mathbb{C}^{n-k}$ then $S_{\Phi}^{-1}\Phi$ is a non-unique $SOD$-frame of $\Phi$ for $r$-erasures.
 \end{thm}
\proof
By Theorem \ref{thm2.8}, the result is true for $r=1,2$, thus, assume $r\geq 3$. Let $G_1,G_2,\ldots,G_k$ be the connected components of $G$ with $|V(G_i)|=n_i$ for $i\in[k]$. By Theorem \ref{thm3.3}, there exist an $L_G(n,n-k)$-frame $\widetilde{\Phi}=\{\widetilde{\Phi}_i\}_{i\in[n]}$ is a non-unique $SOD$-frame for $r$-erasures for any $r\in[n]$ and $\rho_{\widetilde{\Phi},S^{-1}_{\widetilde{\Phi}}\widetilde{\Phi}}^{r}=1$ for $r\in[n-1]\setminus\{1\}$.  

By Theorem \ref{thm2.3}, $\Phi$ and $\widetilde{\Phi}$ are unitary equivalent frames. Hence, there exist a unitary operator $U:\mathbb{C}^{n-k}\rightarrow\mathbb{C}^{n-k}$ such that $U(\widetilde{\phi}_i)=\phi_i$ for all $i\in[n]$. Then, by Theorem \ref{thm3.1}, $\rho^{(r)}_{\Phi,S^{-1}_{\Phi}\Phi}=\rho^{(r)}_{\widetilde{\Phi},S^{-1}_{\widetilde{\Phi}}\widetilde{\Phi}}=1$ for any $r\in[n-1]\setminus\{1\}$.

Let $\Psi=\{\psi_i\}_{i\in[n]}$ be a dual frame of $\Phi$. If $G$ is a null graph then there is nothing to prove. Hence, assume that $G$ is not a null graph. Thus, there exist $j\in[k]$ such that $|V(G_j)|\geq 2$. Let $p,q\in V(G_j) $. Choose $\Lambda\subset[n]$ such that $|\Lambda|=r$ and $p,q\in\Lambda$. By Theorem \ref{thm2.5}, there exist $\nu\in\mathbb{C}^{n-k}$ such that $\psi_i=S_{\Phi}^{-1}\phi_i+\nu$ for $i\in V(G_j)$. Now $E_{\Phi,\Psi,\Lambda}(f)=\sum\limits_{i\in\Lambda}\langle f,\phi_i\rangle \psi_i$. Then, we have
\begin{align*}
   E_{\Phi,\Psi,\Lambda}(\psi_p-\psi_q)&=\sum_{i\in\Lambda} \langle \psi_p-\psi_q,\phi_i\rangle \psi_i\\
   &=\sum_{i\in\Lambda} \langle S_{\Phi}^{-1}\phi_p-S_{\Phi}^{-1}\phi_q,\phi_i\rangle \psi_i\\
    &=\langle S_{\Phi}^{-1}\phi_p-S_{\Phi}^{-1}\phi_q,\phi_i\rangle \psi_p+\langle S_{\Phi}^{-1}\phi_p-S_{\Phi}^{-1}\phi_q,\phi_i\rangle \psi_q\\
    &=\left(1-\frac{1}{n_i}+\frac{1}{n_i}\right)\psi_p+\left(-\frac{1}{n_i}-1+\frac{1}{n_i}\right)\psi_q\\
    &=\psi_p-\psi_q.
\end{align*}
Thus, $1\in\sigma(E_{\Phi,\Psi,\Lambda})$ that is $\rho(E_{\Phi,\Psi,\Lambda})\geq 1$. Hence, $\rho_{\Phi,\Psi}^{(r)}\geq 1=\rho_{\Phi,S^{-1}_{\Phi}\Phi}^{(r)}$.

Since $S_{\Phi}^{-1}\Phi$ is an $SOD$-frame for $\Phi$ for $1$-erasure and $2$-erasures and $\rho_{\Phi,S^{-1}_{\Phi}\Phi}^{(r)}\leq \rho_{\Phi,\Psi}^{(r)}$ for any dual frame $\Psi$ of $\Phi$ and for any $r\in[n-1]$, therefore, $S_{\Phi}^{-1}\Phi$ is an $SOD$-frame for $\Phi$ for $r$-erasures.

Since canonical dual frame of $\widetilde{\Phi}$ is a non-unique optimal dual frame for $r$-erasures, thus, there exist an alternate dual frame $\widetilde\Psi=\{\widetilde{\psi}_i\}_{i\in[n]}$ such that $\widetilde{\Psi}$ is an optimal dual frame for $r$-erasures.  By Theorem \ref{thm2.5}, there exist $\nu_1,\nu_2,\ldots,\nu_k\in\mathbb{C}^{n-k}$ such that 
\begin{align*}
    \widetilde{\Psi}&=\{\widetilde{\psi}_i\}_{i\in[n]}\\
    &=\{S_{\widetilde{\Phi}}^{-1}\widetilde{\phi}_i+\nu_1\}_{i=1}^{n_1}\bigcup\{S_{\widetilde{\Phi}}^{-1}\widetilde{\phi}_i+\nu_2\}_{i=n_1+1}^{n_1+n_2}\bigcup\cdots\bigcup\{S_{\widetilde{\Phi}}^{-1}\widetilde{\phi}_i+\nu_k\}_{i=n_1+\cdots+n_{k-1}+1}^{n_1+\cdots+n_k}
\end{align*}. Choose $\Psi=\{\widetilde{\psi}_i\}_{i\in[n]}=\{S_{\Phi}^{-1}\phi_i+U(\nu_1)\}_{i=1}^{n_1}\bigcup\{S_{\Phi}^{-1}\phi_i+U(\nu_2)\}_{i=n_1+1}^{n_1+n_2}\bigcup\cdots\bigcup\{S_{\Phi}^{-1}\phi_i+U(\nu_k)\}_{i=n_1+\cdots+n_{k-1}+1}^{n_1+\cdots+n_k}$. Then by Theorem \ref{thm2.5}, $\Psi$ is a dual frame of $\Phi$. Also, $\Psi$ is an alternate dual frame of $\Phi$. Further, for any $i\in[n]$, we have
\begin{align*}
S_{\Phi}^{-1}\phi_i+U(\nu_j)&=US_{\widetilde{\Phi}}^{-1}U^*U\widetilde{\phi}_i+U(\nu_j)\\
&=U(S_{\widetilde{\Phi}}^{-1}\widetilde{\phi}_i+\nu_j).
\end{align*}
Hence, $U(\widetilde{\psi}_i)=\psi_i$ for $i\in[n]$. Let $\Lambda\subseteq[n]$ be arbitrary such that $|\Lambda|=r$. Then, for any $f\in\mathbb{C}^{n-k}$, we have
\begin{align*}
    E_{\Phi, \Psi,\Lambda}(f)&=\sum_{i\in\Lambda}\langle f,\phi_i\rangle\psi_i\\
    &=\sum_{i\in\Lambda}\langle f,U(\widetilde{\phi_i})\rangle U(\widetilde{\psi}_i)\\
    &=U(\sum_{i\in\Lambda}\langle U^*(f),\widetilde{\phi_i}\rangle (\widetilde{\psi}_i))\\
    &= UE_{\Phi, \Psi,\Lambda}U^*(f).
\end{align*}
Hence, $ E_{\Phi, \Psi,\Lambda}=U E_{\widetilde{\Phi}, \widetilde{\Psi},\Lambda}U^*$. Thus, we have
\begin{align*}
    \rho_{\Phi, \Psi}^r(f)&=\max\{\rho(E_{\Phi,\Psi,\Lambda}):|\Lambda|=r, \Lambda\subseteq[n]\}\\
    &=\max\{\rho(UE_{\widetilde{\Phi},\widetilde{\Psi},\Lambda}U^*):|\Lambda|=r, \Lambda\subseteq[n]\}\\
    &=\max\{\rho(E_{\widetilde{\Phi},\widetilde{\Psi},\Lambda}):|\Lambda|=r, \Lambda\subseteq[n]\}\\
    &=\rho_{\widetilde{\Phi}, S_{\widetilde{\Phi}}}\widetilde{\Phi}^r
\end{align*}.
Thus, $\Psi$ is an $SOD$-frame of $\Phi$ for $r$-erasures.

\endproof

    \section*{Acknowledgments}
	Aniruddha Samanta expresses thanks to the National Board for Higher Mathematics (NBHM), Department of Atomic Energy, India, for providing financial support in the form of an NBHM Post-doctoral Fellowship (Sanction Order No. 0204/21/2023/R\&D-II/10038). The second author also acknowledges excellent working conditions in the Theoretical Statistics and Mathematics Unit, Indian Statistical Institute Kolkata. 
	
	\mbox{}

\end{document}